\documentclass[12pt]{article}
\usepackage[T2A]{fontenc}
\usepackage{amsmath}
\usepackage[cp1251]{inputenc}
\begin{document}

\vspace{0.2cm}

\large

\begin{center}
{\bf CAN ONE HEAR FASTENING OF A ROD?}
\end{center}

\centerline{\bf Akhtyamov A.M., Mouftakhov A.V.,}

\centerline{\bf Teicher M., Yamilova L.S.}

\centerline{\bf (Bashkir State University, Ufa;}

\centerline{\bf Bar-Ilan University, Ramat-Gan,Israel)}

\vspace{0.2cm}

\vspace{0.2cm}

{\bf 1. Introduction.} Rods are parts of various devices (see
\cite{Timoshenko 55}-\cite{Vibrations 78}). If it is impossible to
observe the rod directly, the only source of information about
possible defects of its fastening can be the natural frequencies of
its flexural vibrations. The question arises whether one would be
able to detect damage in rod fastening by the natural frequencies of
its flexural vibrations. This paper gives and substantiates a
positive answer to this question.

 The problem in question
belongs to the class of inverse problems and is a completely
natural problem of identification of the boundary conditions.

Closely related formulations of the problem were proposed in
\cite{Kac 66,Qunli 90}. Contrary to this, in this paper it is not
the form of the domain or size of an object which are sought for
but the nature of fastening. The problem of determining a boundary
condition has been considered in \cite{Frikha 00}. However, as
data for finding the boundary conditions, we take not a set of
natural frequencies, but not condensation and inversion (as in
\cite{Frikha 00}).

Similarly formulated problems also occur in the spectral theory of
differential operators, where it is required to determine the
coefficients of a differential equation and the boundary
conditions using a set of eigenvalues (for more details, see
\cite{Borg 46}--\cite{Akhtyamov 99 DE}). However, as data for
finding the boundary conditions, we take  one spectrum  but not
several spectra or other additional spectral data (for example,
the spectral function, the Weyl function or the so-called
weighting numbers) that were used in \cite{Borg
46}--\cite{Akhtyamov 04 RAN a}. Moreover, the principal aim there
was to determine the coefficients in the equation and not in
boundary conditions. The aim of this paper is to determine  the
boundary conditions of the eigenvalue problem from its spectrum in
the case of a known differential equation.

The problem of determining a boundary condition using a finite set
of eigenvalues has been considered previously in
  \cite{Akhtyamov 03 AJ-a}--\cite{Akhtyamov 01}. In contrast to
  papers \cite{Akhtyamov 03 AJ-a}--\cite{Akhtyamov 03 IPE},
  we think it is necessary to determine not
the  type of fastening of plates, but the type of fastening of a
rod. The problem of finding of fastening of one of the end of a
rod is considered in \cite{Akhtyamov 01}. In contrast to
  papers \cite{Akhtyamov 01} we think it is necessary to determine not
the  type of fastening of one of the ends of a rod, but the  type
of fastening of both ends of a rod.

 \vspace{0.2cm}

{\bf 2. Formulation of the direct problem.}  The problem of
flexural oscillations of rod is reduced to following boundary
problem \cite{Timoshenko 55,Collatz 68 a}
\begin{equation}\label{rod1 1}
(\alpha\, u'')''=\rho\, F\, \ddot{u},
\end{equation}
\begin{equation}\label{rod1 2}
\left[ (\alpha\, u'')'-c_1^-\, u\right]_{x=0}=0, \qquad \left[
\alpha\, u''-c_2^-\, u\right]_{x=0}=0,
\end{equation}
\begin{equation}\label{rod1 3}
\left[ (\alpha\, u'')'+c_1^+\, u'\right]_{x=l}=0, \qquad
\left[
\alpha\, u''+c_2^+\, u'\right]_{x=l}=0,
\end{equation}
 where $\alpha$ is the flexural rigidity, $l$ is the length, $\rho$ is the density
and $F$ is the cross-section area of the rod, $c_i^\pm$ ($0\leq
c_i^\pm\leq \infty$) are rigidity coefficients of springs.

Let $\alpha$, $\rho$ and $F$ be constant. By $\lambda^2$ denote
 $\rho\, F\, \omega^2\, /\, \alpha$.

For vibrations, we write
$$
u(x,t)=y(x)\, \cos {\omega\, t}
$$
and hence obtain the following eigenvalue problem for $y(x)$:
\begin{equation}\label{rod1 1'}
y^{(4)}=s^2 y,
\end{equation}
\begin{equation}\label{rod1 2'}
U_1(y)=-a_1\, y(0)+a_4\, y'''(0)=0, \quad U_2(y)=-a_2\,
y'(0)+a_3\, y''(0)=0,
\end{equation}
\begin{equation}\label{rod1 3'}
U_3(y)=a_5\, y(l)+a_8\, y'''(l)=0,\quad U_4(y)=a_6\, y'(l)+a_7\,
y''(l)=0,
\end{equation}
where $s^2=\rho\, F\, \omega^2/\alpha$, \ $a_i\geq 0$ \
($i=1,2,\dots,8$).

The coefficients $a_i$ characterize conditions for fastening the
rod (rigid clamping, free support, free edge, floating fixing,
elastic fixing).

We note that the functions
\begin{equation}\label{rod1 y1234}
\begin{array}{l}
  y_1(x,s )= (\cos\sqrt{s}\,  x +\cosh \sqrt{s}\,  x)/2,
  \\ y_2(x,s )= (\sin\sqrt{s}\,  x +\sinh \sqrt{s}\,  x)/(2\sqrt{s} ),
  \\
y_3(x,s )= (-\cos\sqrt{s}\,  x +\cosh \sqrt{s}\, x)/(2s),
\\
y_4(x,s )= (-\sin\sqrt{s}\,  x +\sinh \sqrt{s}\, x)/(2\sqrt{s^3}),
\end{array}
\end{equation}
are linearly independent solutions of the equation (\ref{rod1 1'})
which satisfy the conditions
\begin{equation}\label{rod1 delta1}
  y_j^{(r-1)}(0,s)=\left\{
  \begin{array}{l}
0\quad \mbox{when}\quad j\ne r,
\\
1\quad \mbox{when}\quad j=r,
\end{array}
\qquad j,r=1,2,3,4. \right.
\end{equation}

The natural frequencies $\omega_i$ ($s_i$) are the corresponding
positive eigenvalues of problem (\ref{rod1 1'})--(\ref{rod1 3'})
(see \cite{Collatz 68 a,Vibrations 78}). The non-zero eigenvalues
of problem (\ref{rod1 1'})--(\ref{rod1 3'}) are the roots of the
determinant
$$
   \Delta (s)=
\left|\begin{array}{cccc}
  U_1(y_1) & U_1(y_2) & U_1(y_3) & U_1(y_4) \\
  U_2(y_1) & U_2(y_2) & U_2(y_3) & U_2(y_4) \\
  U_3(y_1) & U_3(y_2) & U_3(y_3) & U_3(y_4) \\
  U_4(y_1) & U_4(y_2) & U_4(y_3) & U_4(y_4)
\end{array}\right|
$$

We shall now formulate the direct eigenvalue problem (\ref{rod1
1'})--(\ref{rod1 3'}): it is required to find the unknown natural
frequencies of the oscillations of the rod from
 $a_i$ ($i=1,2,\dots 8$).

Thus, finding the natural frequencies $\omega_i$ ($i=1,2,\dots $)
is equivalent to founding of the roots $s_i$ ($i=1,2,\dots $) of
$\Delta (s )$.

Thus, knowing $a_i$ ($i=1,2,\dots 8$) it is possible to find $s_i$
($i=1,2,\dots $) by standard methods \cite{Collatz 68 a,Vibrations
78}. The solution to this direct problem presents no difficulties.
The question arises whether one would be able to do the reverse
and find $a_i$ ($i=1,2,\dots 8$) knowing $s_i$ ($i=1,2,\dots $).
In a broader sense it may be stated as follows. Is it possible to
derive unknown boundary conditions with a knowledge of $s_i$? The
answer to this question is given in the next section.

\vspace{0.2cm}

{\bf 3. Formulation of the inverse problem.} We shall denote the
matrix, consisting of the coefficients $a_{j}$ ($j=1,2,\dots ,8$)
of the forms $U_i(y)$ ($i=1,2,3,4$) by $A$ and its minors by
$M_{ijkm}$ ($i,j,m,k=1,2,\dots ,8$):
$$
A= \left\|\begin{array}{cccccccc} a_{1} & 0 & 0 & a_4 & 0 & 0 & 0
& 0
\\
0 & a_{2} & a_{3} & 0 & 0 & 0 & 0 & 0
\\
0 & 0 & 0 & 0 & a_{5} & 0 & 0 & a_{8}
\\
0 & 0 & 0 & 0 & 0 & a_{6} & a_{7} & 0
\end{array} \right\| ,
\qquad M_{ijkm}=\pm\, a_{i}\, a_{j}\, a_{k} \, a_{m}.
$$

In terms of eigenvalue problem (\ref{rod1 1'})--(\ref{rod1 3'}),
the inverse problem constructed above should be formulated as
follows: the coefficients $a_{j}$ ($j=1,2,\dots 8$) of the forms
$U_i(y)$ ($i=1,2,3,4$)  of problem (\ref{rod1 1'})--(\ref{rod1
3'}) are unknown, the rank of the matrix A to make up these
coefficients is equal to four, the eigenvalues $s_k$ of problem
(\ref{rod1 1'})--(\ref{rod1 3'}) are known and it is required to
find the matrix $A$ (the class of linear equivalent matrixes).

\vspace{0.2cm}

{\bf 4. The duality
 of the solution.}
 Together with problem (\ref{rod1 1'})--(\ref{rod1 3'}), let us
consider the following eigenvalue problem
\begin{equation}\label{rod1 1''}
y^{(4)}=s^2 y,
\end{equation}
\begin{equation}\label{rod1 2''}
\widetilde{U_1}(y)=-\widetilde{a_1}\, y(0)+\widetilde{a_4}\,
y'''(0)=0, \quad \widetilde{U_2}(y)=-\widetilde{a_2}\,
y'(0)+\widetilde{a_3}\, y''(0)=0,
\end{equation}
\begin{equation}\label{rod1 3''}
\widetilde{U_3}(y)=\widetilde{a_5}\, y(l)+\widetilde{a_8}\,
y'''(l)=0,\quad \widetilde{U_4}(y)=\widetilde{a_6}\,
y'(l)+\widetilde{a_7}\, y''(l)=0.
\end{equation}

We denote the matrix composed of the coefficients
$\widetilde{a_{i}}$ of the forms $\widetilde{U_j}(y)$ by
$\widetilde{A}$ and its minors by $\widetilde M_{ijkm}$:
$$
\widetilde{A}= \left\|\begin{array}{cccccccc} \widetilde{a_{1}} &
0 & 0 & \widetilde{a_4} & 0 & 0 & 0 & 0
\\
0 & \widetilde{a_{2}} & \widetilde{a_{3}} & 0 & 0 & 0 & 0 & 0
\\
0 & 0 & 0 & 0 & \widetilde{a_{5}} & 0 & 0 & \widetilde{a_{8}}
\\
0 & 0 & 0 & 0 & 0 & \widetilde{a_{6}} & \widetilde{a_{7}} & 0
\end{array} \right\| ,
\qquad \qquad \widetilde{M}_{ijkm}=\pm\, \widetilde{a_{i}}\,
\widetilde{a_{j}}\, \widetilde{a_{k}} \, \widetilde{a_{m}}.
$$

Let $\widetilde A^*$ be matrix
$$
\widetilde A^*=  \left\|\begin{array}{cccccccc} \widetilde{a_{5}}
& 0 & 0 & \widetilde{a_8} & 0 & 0 & 0 & 0
\\
0 & \widetilde{a_{6}} & \widetilde{a_{7}} & 0 & 0 & 0 & 0 & 0
\\
0 & 0 & 0 & 0 & \widetilde{a_{1}} & 0 & 0 & \widetilde{a_{4}}
\\
0 & 0 & 0 & 0 & 0 & \widetilde{a_{2}} & \widetilde{a_{3}} & 0
\end{array} \right\| .
$$

{\it Theorem 1 (on the duality of the solution of the inverse
problem).} Suppose the following conditions are satisfied
\begin{equation}\label{ipi 2 1}
  {\rm rank\, } A={\rm rank\, } \widetilde A =4.
\end{equation}
If the eigenvalues $\{ s_i\}$ of problem (\ref{rod1
1'})--(\ref{rod1 3'}) and the  eigenvalues $\{ \widetilde s_i\}$
of problem (\ref{rod1 1''})--(\ref{rod1 3''}) are identical, with
account taken for their multiplicities, the classes of linearly
equivalent matrixes $A$ and $\widetilde A$ or $A$ and $\widetilde
A^*$  are also identical.

{\it Proof.}
 The eigenvalues of problem (\ref{rod1
1'})--(\ref{rod1 3'}) are the roots of the determinant

\begin{equation}\label{rod1 Delta1}
\begin{array}{l}  \Delta (s)=
  -M_{1256}\,[f^-(s)/s^2]+(M_{2457}+M_{1368})\, [f^-(s)]+M_{3478}\, [s^2\,f^-(s)]
  \\+(M_{1278}+  M_{3456}-M_{2457}-M_{1368})\, [f^+(s)]+M_{1357}\,
  [z(s)/s]-M_{2468}\, [s\, z(s)]\\
+(M_{2456}+M_{1268})\,[g^-(s)/\sqrt{s}]-(M_{3468}+M_{2478})\,
    [\sqrt{s^3}\,g^-(s)]+
    \\
    +(M_{1356}+M_{1257})\,[ g^+(s)/\sqrt{s^3}]-(M_{3457}+M_{1378})\, [\sqrt{s}\, g^+(s)],
\end{array}
\end{equation}
where
$$
f^\pm(s)=(1\pm\cos\sqrt{s}\, \cosh \sqrt{s})/2, \quad z(s)=
(\sin\sqrt{s}\,\sinh\sqrt{s})/2,
$$
$$g^\pm(s)= (-\sin \sqrt{s}\, \cosh\sqrt{s} \pm \cos\sqrt{s}\,
\sinh\sqrt{s})/2.
$$
Note that the functions in square brackets are linearly
independent. Since $\Delta (s)\not\equiv 0$, $\widetilde\Delta
(s)\not\equiv 0$ are entire functions in s of order 1/2, it
follows from Hadamard's factorization theorem (see \cite{Gol'dberg
91}) that determinants $\Delta (s)$  and $\widetilde\Delta (s)$
are connected by the relation
$$
\Delta (s)\equiv C\,  \widetilde\Delta (s),
$$
where $k$ is a certain  integer and $C$ is a certain non-zero
constant.  From this and we obtain the equalities
\begin{gather}
M_{1256}=K\, \widetilde M_{1256}, \label{rod1 M25}
\\
M_{1357}=K\, \widetilde M_{1357}, \label{rod1 M26}
\\
M_{2468}=K\, \widetilde
 M_{2468},\label{rod1 M27}
\\
M_{3478}=K\, \widetilde M_{3478},\label{rod1 M28}
\\
  M_{1257}+M_{1356}=K\, (\widetilde M_{1356}+\widetilde
M_{1257}),\label{rod1 M29}
\\
M_{1268}+M_{2456}=K\,(\widetilde M_{1268}+\widetilde
M_{2456}),\label{rod1 M30}
\\
M_{1378}+M_{3457}=K\, (\widetilde M_{1378}+\widetilde
M_{3457}),\label{rod1 M31}
\\
M_{2478}+M_{3468}=K\, (\widetilde M_{2478}+\widetilde
M_{3468}),\label{rod1 M32}
\\
M_{1278}+M_{3456}=K\, (\widetilde M_{1278}+\widetilde
M_{3456}),\label{rod1 M33}
\\
M_{1368}+M_{2457}=K\, \,(\widetilde M_{1368}+\widetilde
M_{2457}).\label{rod1 M34}
\end{gather}

Hence, \[
a_1=\widetilde{a_1}, \quad a_2=\widetilde{a_2}, \quad
a_3=\widetilde{a_3}, \quad a_4=\widetilde{a_4},
\]
\[ a_5=\widetilde{a_5}, \quad a_6=\widetilde{a_6},  \quad
a_7=\widetilde{a_7},  \quad a_8=\widetilde{a_8},
\]
or
\[
a_1=\widetilde{a_5}, \quad a_2=\widetilde{a_6}, \quad
a_3=\widetilde{a_7}, \quad a_4=\widetilde{a_8},
\]
\[ a_5=\widetilde{a_1}, \quad a_6=\widetilde{a_4},  \quad
a_7=\widetilde{a_3},  \quad a_8=\widetilde{a_4}.
\]

These equations can be proved by standard method of algebra. Let
us consider only one case $M_{3478}\ne 0$ for example.  If
$M_{3478}\ne 0$,
 $a_3\ne 0$, $a_4\ne 0$, $a_7\ne 0 $, $a_8\ne 0$,
$\widetilde{a_3}\ne 0$, $\widetilde{a_4}\ne 0$,
$\widetilde{a_7}\ne 0 $, $\widetilde{a_8}\ne 0$ and so matrixes
$A$ and $\widetilde{A}$ have the forms
\[\label{rod1 CAB} A
 = \left\|
\begin{array}{cccccccc}
a_1&0&0&1&0&0&0&0
\\
0&a_2&1&0&0&0&0&0
\\
0&0&0&0&a_5&0&0&1
\\
0&0&0&0&0&a_6&1&0
\end{array}
\right\| ,\quad \widetilde{A}
 = \left\|
\begin{array}{cccccccc}
\widetilde{a_1}&0&0&1&0&0&0&0
\\
0&\widetilde{a_2}&1&0&0&0&0&0
\\
0&0&0&0&\widetilde{a_5}&0&0&1
\\
0&0&0&0&0&\widetilde{a_6}&1&0
\end{array}
\right\| . \]

From this and from (\ref{rod1 M28}) we obtain $K=1.$

It follows from (\ref{rod1 M26}), (\ref{rod1 M31}) that
$$a_1\, a_5=\widetilde{a_1}\, \widetilde{a_5},
\quad  a_1+a_5=\widetilde{a_1}+\widetilde{a_5}.$$ Therefore
 by Vieta's theorem, we have
 \begin{equation}\label{prpl-ar5 35}
 a_1=\widetilde{a_1},\quad a_5=\widetilde{a_5}
 \quad \mbox{ or }\quad
a_1=\widetilde{a_5},\quad  a_5=\widetilde{a_1}.
\end{equation}

Similarly, equations (\ref{rod1 M27}), (\ref{rod1 M32}) imply
\begin{equation}\label{rod1 36}
 a_2=\widetilde{a_2},\quad a_6=\widetilde{a_6},
 \quad \mbox{ or }\quad
a_2=\widetilde{a_6},\quad  a_6=\widetilde{a_2}.
\end{equation}

These equations are equivalent to four cases:

Case 1:  $a_1=\widetilde{a_1},\quad a_5=\widetilde{a_5},\quad
 a_2=\widetilde{a_2},\quad a_6=\widetilde{a_6}$.

Case 2: $a_1=\widetilde{a_5},\quad a_5=\widetilde{a_1},\quad
a_2=\widetilde{a_6},\quad  a_6=\widetilde{a_2}$.

Case 3: $a_1=\widetilde{a_1},\quad a_5=\widetilde{a_5},\quad
a_2=\widetilde{a_6},\quad  a_6=\widetilde{a_2}$.

Case 4:  $a_1=\widetilde{a_5},\quad  a_5=\widetilde{a_1},\quad
  a_2=\widetilde{a_2},\quad a_6=\widetilde{a_6}$.

Cases 3 and 4 are special cases of 1 and 2. Indeed,  in case 3 it
follows from (\ref{rod1 M33}) or (\ref{rod1 M34}) that
\begin{equation}\label{rod1 37}
  (a_1-a_5)(a_2-a_6)=0.
\end{equation}
Hence $a_1-a_5=0$ or $a_2-a_6=0.$ If $a_1-a_5=0$, then case 3 is
special case 2: $\widetilde{a_1}=a_1=a_5=\widetilde{a}_5$. If
$a_2-a_6=0$, then case 3 is special case 1:
$\widetilde{a_6}=a_2=a_6=\widetilde{a}_6$.

In the same way, we can show that
  the case 4 is special case of cases 1 and 2. We have proved that
in the case $M_{3478}\ne 0$ it is possible only two cases: 1 and
2.  In a similar manner we can prove the analogous results for
other cases. This contradiction proves the theorem.

\vspace{0.2cm}

{\bf 5. Exact solution.} It has been shown above that the problem
of finding the matrix $A$ (in the sense that class of linearly
equivalent matrixes of $A$) from the natural frequencies of
flexural oscillations of a rod has a duality solution. The next
question is how this solution can be constructed.

This section deals with solving this problem and constructing
exact solution by the first 9 natural frequencies $\omega _i$.

Suppose  $s_1$,   $s_2$, \dots , $s_9$ are the values
corresponding to the first nine natural frequencies $\omega _i$.
We substitute the values $s_1$,   $s_2$, \dots , $s_9$ into
(\ref{rod1 Delta1}) and obtain a system of nine homogeneous
algebraic equations
\begin{equation}\label{rod1 Delta2}
\begin{array}{l}
  x_1\,[f^-(s_i)/s_i^2]+x_2\, [f^-(s_i)]+x_3\, [s_i^2\,f^-(s_i)]
  \\+x_4\, [f^+(s_i)]+x_5\,
  [z(s_i)/s_i]+x_6\, [s_i\, z(s_i)]\\
+x_7\,[g^-(s_i)/\sqrt{s_i}]+x_8\,
    [\sqrt{s_i^3}\,g^-(s_i)]+
    \\
    +x_{9}\, g^+(s_i)/\sqrt{s_i^3}]+x_{10}\, [\sqrt{s_i}\,
    g^+(s_i)]=0,\qquad i=1,2,\dots , 9
\end{array}
\end{equation}
 in the ten unknowns
\begin{equation}\label{rod1 x}
\begin{array}{l}
x_1=
  -M_{1256},\quad x_2=M_{2457}+M_{1368}, \quad x_3=M_{3478},
  \\
 x_4=M_{1278}+  M_{3456}-M_{2457}-M_{1368},\quad
 x_5=M_{1357},\\
   x_6=-M_{2468},\quad
  x_7=
M_{2456}+M_{1268},\quad  x_8=-(M_{3468}+M_{2478}),\\
x_9=M_{1356}+M_{1257}, \quad  x_{10}=-(M_{3457}+M_{1378}).
\end{array}
\end{equation}

 The
resulting set of equations (\ref{rod1 Delta2})  has an infinite
number of solutions. If the resulting set has a rank of 9, the
unknown $x_i$ can be found in accurate to a coefficient.  The
unknown matrix $A$ is found from (\ref{rod1 x}) by direct
calculations as in proof of theorem~1.

These reasons prove

{\it Theorem 2 (on the duality of the solution of the inverse
problem).} {If the matrix of system (\ref{rod1 Delta2}) has a rank
of 9, the solution of the inverse problem of the reconstruction
boundary conditions~(\ref{rod1 2''}), (\ref{rod1 3''}) is
duality.}

{\it Remark.} Theorem 2 is stronger than theorem 1. Theorem 2 use
only 9 natural frequencies for the reconstruction of boundary
conditions and
 not all natural frequencies as in
theorem 1.
   But theorem 1 is proves duality of the solution in the common case
   (if the set (\ref{rod1 Delta2}) has not a rank of 9).

 Continuity of the solution
 of the inverse problem with respect to $s_i$ is proved as in \cite{Akhtyamov 03 IPE}.
 This shows that small perturbations of eigenvalues
 $s_i$ ($i=1,2,3$) lead to small perturbations of the
 boundary conditions.
It follows from this and theorem~1  that the inverse problem is
well posed, since its solution exists, is unique and
 continuous with respect to $s_i$ $(i=1,2,3)$.

 Computer calculations confirm the stability of the solution
 of the inverse problem. The order of error often
hardly different from the error in the closeness of values of
$\widetilde s_k$ and $s_i$ and only in some cases it can be
deteriorated by four orders of magnitude. So the measurement
accuracy of instruments to measure natural frequencies must exceed
accuracy to measure boundary conditions by four orders of
magnitude.

\vspace{0.2cm}

{\bf 3. Numerical results.}  We use dimensionless variables in the
numerical examples.

\vspace{0.2cm}

 {\it Example 1 (rigid clamping --- free support).} Suppose
$$
\begin{array}{c}
s_1=15.4182057169801,\quad s_2=49.9648620318002, \\
s_3=104.247696458861, \quad
 s_4=178.269729494609,\\ s_5=272.030971305025, \quad s_6=385.531421917553,
 \\
  s_7=518.771081332259, \quad s_8=671.749949549144, \\
   s_9=844.468026568208
\end{array}
$$
 correspond to the first 9 natural
frequencies $\omega_i$
 determined using instruments for measuring the natural frequencies,
 then the solution of set
 (\ref{rod1 Delta2}) has the form
\begin{equation}\label{rod1 1prx}
\begin{array}{l}
x_1=
  -M_{1256}\approx 0\cdot C,\quad x_2=M_{2457}+M_{1368}\approx 0\cdot C, \\
x_3=M_{3478}\approx 0\cdot C,
  \quad
 x_4=M_{1278}+  M_{3456}-M_{2457}-M_{1368}\approx 0\cdot C,\\
 x_5=M_{1357}\approx 0\cdot C,\quad
   x_6=-M_{2468}\approx 0\cdot C,\\
  x_7=-M_{2456}+M_{1268}\approx 0\cdot C,\quad
x_8=-(M_{3468}+M_{2478})\approx 0\cdot C,
\\
x_9=M_{1356}+M_{1257}\approx  C , \quad
x_{10}=-(M_{3457}+M_{1378})\approx 0\cdot C
\end{array}
\end{equation}
 with an accuracy of $10^{-9}$.

Substituting  $1$ for $C$ in (\ref{rod1 1prx}), we get
$M_{1356}+M_{1257}=1$. This means that $M_{1356}\ne 0$ or
$M_{1257}\ne 0$.

Suppose $M_{1356}\ne 0$. Then $a_1\, a_3\, a_5\, a_6\ne 0$. So
$$
A= \left\|\begin{array}{cccccccc} 1 & 0 & 0 & a_4 & 0 & 0 & 0 & 0
\\
0 & a_2 & 1 & 0 & 0 & 0 & 0 & 0
\\ 0 & 0 & 0 & 0 & 1 & 0 & 0 & a_8
\\ 0 & 0 & 0 & 0 & 0 & 1 & a_7 & 0
\end{array}\right\|
$$
(in the sense that class of linearly equivalent matrixes of $A$).

Using $M_{1256}=0$, $M_{1357}=0$ and $M_{2457}+M_{1368}=0$, we
have $a_2=0$, $a_7=0$ and $a_8=0$. From this and
$M_{1278}+M_{3456}-M_{2457}-M_{1368}=0$, we get $a_4=0$.

Thus we have
$$
A= \left\|\begin{array}{cccccccc} 1 & 0 & 0 & 0 & 0 & 0 & 0 & 0
\\
0 & 0 & 1 & 0 & 0 & 0 & 0 & 0
\\ 0 & 0 & 0 & 0 & 1 & 0 & 0 & 0
\\ 0 & 0 & 0 & 0 & 0 & 1 & 0 & 0
\end{array}\right\|
$$

Suppose now $M_{1257}\ne 0$.
 Then $a_1\, a_2\, a_5\, a_7\ne 0$.
 So
$$
A= \left\|\begin{array}{cccccccc} 1 & 0 & 0 & a_4 & 0 & 0 & 0 & 0
\\
0 & 1 & a_3 & 0 & 0 & 0 & 0 & 0
\\ 0 & 0 & 0 & 0 & 1 & 0 & 0 & a_8
\\ 0 & 0 & 0 & 0 & 0 & a_6 & 1 & 0
\end{array}\right\|
$$

Using $M_{1256}=0$ and $M_{2457}+M_{1368}=0$, we get $a_6=0$ and
$a_4=0$.

If we combine this with equations $M_{1357}=0$ and
 $M_{1278}+M_{3456}-M_{2457}-M_{1368}=0$,  we get $a_3=0$ and $a_8=0$.

Thus we have
$$
A= \left\|\begin{array}{cccccccc} 1 & 0 & 0 & 0 & 0 & 0 & 0 & 0
\\
0 & 1 & 0 & 0 & 0 & 0 & 0 & 0
\\ 0 & 0 & 0 & 0 & 1 & 0 & 0 & 0
\\ 0 & 0 & 0 & 0 & 0 & 0 & 1 & 0
\end{array}\right\|
$$

 So
boundary conditions are or
$$
U_1(y)=y(0)=0, \quad U_2(y)= y''(0)=0,
$$
$$
U_3(y)=y(l)=0,\quad U_4(y)=y'(l)=0,
$$
 either
$$
U_1(y)= y(0)=0,\quad U_2(y)= y'(0)=0,
$$
$$
U_3(y)=y(l)=0, \quad U_4(y)= y''(l)=0.
$$

 Note that the  numbers ${s_i}$
 presented above are the same as the first
nine exact values corresponding to fastening of rigid clamping and
free support. This means that the unknown rod fastening
inaccessible to direct observation has been correctly determined.

\vspace{0.2cm}

 {\it Example 2 (rigid clamping --- elastic
fixing).} If
$$
\begin{array}{c}
s_1=5.60163863016235,\quad s_2=22.4984332740862, \\
s_3=61.8604321649037, \quad
 s_4=120.984868139371,\\ s_5=199.909638169628, \quad s_6=298.589053349029,
 \\
  s_7=417.014779762035, \quad s_8=555.183266366176, \\
   s_9=713.092945010199
\end{array}
$$
 correspond to the first 9 natural
frequencies $\omega_i$
 determined using instruments for measuring the natural frequencies,
 then the solution of set
 (\ref{rod1 Delta2}) has the form
\begin{equation}\label{rod1 x}
\begin{array}{l}
x_1=
  -M_{1256}\approx 0,
  \quad x_2=M_{2457}+M_{1368}\approx 0, \\
x_3=M_{3478}\approx 0,
  \quad
 x_4=M_{1278}+  M_{3456}-M_{2457}-M_{1368}\approx -C,\\
 x_5=M_{1357}\approx 0,\quad
   x_6=-M_{2468}\approx 0,\\
  x_7=-M_{2456}+M_{1268}\approx 0,\quad
x_8=-(M_{3468}+M_{2478})\approx 0
,\\
x_9=M_{1356}+M_{1257}\approx 5\, C , \quad
x_{10}=-(M_{3457}+M_{1378})\approx 0.
\end{array}
\end{equation}
Suppose  $C=1$; then $M_{1356}+M_{1257}\ne 0$. This means that
$M_{1356}\ne 0$ or $M_{1257}\ne 0$.

Suppose $M_{1356}\ne 0$. Then $a_1\, a_3\, a_5\, a_6\ne 0$. So
$$
A= \left\|\begin{array}{cccccccc} 1 & 0 & 0 & a_4 & 0 & 0 & 0 & 0
\\
0 & a_2 & 1 & 0 & 0 & 0 & 0 & 0
\\ 0 & 0 & 0 & 0 & 1 & 0 & 0 & a_8
\\ 0 & 0 & 0 & 0 & 0 & 1 & a_7 & 0
\end{array}\right\|
$$
(in the sense that class of linearly equivalent matrixes of $A$).

If we combine this with equations $M_{1256}=0$,
$M_{1356}+M_{1257}=5\, C$ we get $a_2=0$, $1=5\, C$. Using
$M_{1357}=0$ and $M_{2457}+M_{1368}=0$, we have $a_7=0$ and
$a_8=0$. From this and $M_{1278}+M_{3456}-M_{2457}-M_{1368}=-C$,
we get $-a_4=-1/5$.

Thus we have
$$
A= \left\|\begin{array}{cccccccc} 1 & 0 & 0 & 1/5 & 0 & 0 & 0 & 0
\\
0 & 0 & 1 & 0 & 0 & 0 & 0 & 0
\\ 0 & 0 & 0 & 0 & 1 & 0 & 0 & 0
\\ 0 & 0 & 0 & 0 & 0 & 1 & 0 & 0
\end{array}\right\|
$$

Suppose now $M_{1257}\ne 0$.
 Then $a_1\, a_2\, a_5\, a_7\ne 0$.
 So
$$
A= \left\|\begin{array}{cccccccc} 1 & 0 & 0 & a_4 & 0 & 0 & 0 & 0
\\
0 & 1 & a_3 & 0 & 0 & 0 & 0 & 0
\\ 0 & 0 & 0 & 0 & 1 & 0 & 0 & a_8
\\ 0 & 0 & 0 & 0 & 0 & a_6 & 1 & 0
\end{array}\right\|
$$

Using $M_{1256}=0$ and $M_{2457}+M_{1368}=0$, we get $a_6=0$ and
$a_4=0$.

If we combine this with equations $M_{1357}=0$,
$M_{1356}+M_{1257}=5\, C$,
 $M_{1278}+M_{3456}-M_{2457}-M_{1368}=-C$ and
$M_{2457}+M_{1368}=0$ we get $a_3=0$, $1=5\, C$, and $-a_8=-1/5$.

Thus we have
$$
A= \left\|\begin{array}{cccccccc} 1 & 0 & 0 & 0 & 0 & 0 & 0 & 0
\\
0 & 1 & 0 & 0 & 0 & 0 & 0 & 0
\\ 0 & 0 & 0 & 0 & 1 & 0 & 0 & 1/5
\\ 0 & 0 & 0 & 0 & 0 & 0 & 1 & 0
\end{array}\right\|
$$

 So
boundary conditions are or
$$
U_1(y)=- 5\, y(0)+y'''(0)=0, \quad U_2(y)= y''(0)=0,
$$
$$
U_3(y)=y(l)=0,\quad U_4(y)=y'(l)=0,
$$
 either
$$
U_1(y)= y(0)=0,\quad U_2(y)= y'(0)=0,
$$
$$
U_3(y)=5\, y(l)+y'''(l)=0, \quad U_4(y)= y''(l)=0.
$$

 Note that the  numbers ${s_i}$
 presented above are the same as the first
nine exact values corresponding to fastening of rigid clamping
and elastic fixing with the relative stiffness factor of 5. This
means that the unknown rod fastening inaccessible to direct
observation has been correctly determined.

\vspace{0.2cm}

{\it Example 3. (elastic fixing --- elastic fixing)} If
$$
\begin{array}{c}
s_1=0.383848559322840,\quad s_2=11.9180148367849, \\
s_3=53.4326824121208, \quad
 s_4=114.157790468867,\\ s_5=193.836586296759, \quad s_6=292.955617117120,
 \\
  s_7=411.667770695782, \quad s_8=550.037492353361, \\
   s_9=708.096219400352
\end{array}
$$
 correspond to the first 9 natural
frequencies $\omega_i$
 determined using instruments for measuring the natural frequencies
 with an accuracy of $10^{-15}$,
 then the solution of set
 (\ref{rod1 Delta2}) has the form
\begin{equation}\label{rod1 x}
\begin{array}{l}
x_1=
  -M_{1256}\approx -24\, C
,\quad x_2=M_{2457}+M_{1368}\approx -10\, C , \\
x_3=M_{3478}\approx C,
  \quad
 x_4=M_{1278}+  M_{3456}-M_{2457}-M_{1368}\approx 9\, C
,\\
 x_5=M_{1357}\approx 3\, C
,\quad
   x_6=-M_{2468}\approx -8\, C
,\\
  x_7=-M_{2456}+M_{1268}\approx -32\, C ,\quad
x_8=-(M_{3468}+M_{2478})\approx -6\, C
,\\
x_9=M_{1356}+M_{1257}\approx 18\, C , \quad
x_{10}=-(M_{3457}+M_{1378})\approx 4\, C.
\end{array}
\end{equation}
Suppose  $C=1$; then $M_{3478}=a_3\, a_4\, a_7\, a_8\ne 0$. So
$$
A= \left\|\begin{array}{cccccccc} a_1 & 0 & 0 & 1 & 0 & 0 & 0 & 0
\\
0 & a_2 & 1 & 0 & 0 & 0 & 0 & 0
\\ 0 & 0 & 0 & 0 & a_5 & 0 & 0 & 1
\\ 0 & 0 & 0 & 0 & 0 & a_6 & 1 & 0
\end{array}\right\|
$$
(in the sense that class of linearly equivalent matrixes of $A$).

From this and (\ref{rod1 x}), we obtain or $a_1=1$, $a_2=2$,
$a_5=3$, $a_6=4$ either $a_1=3$, $a_2=4$, $a_5=1$, $a_6=2$. So
boundary conditions are or
$$
U_1(y)=-y(0)+y'''(0)=0, \quad U_2(y)=-2\, y'(0)+ y''(0)=0,
$$
$$
U_3(y)=3\, y(l)+ y'''(l)=0,\quad U_4(y)=4\, y'(l)+y''(l)=0,
$$
 either
$$
U_1(y)=-3\, y(0)+ y'''(0)=0,\quad U_2(y)=-4\, y'(0)+y''(0)=0,
$$
$$
U_3(y)=y(l)+y'''(l)=0, \quad U_4(y)=2\, y'(l)+ y''(l)=0.
$$

 Note that the  numbers ${s_i}$
 presented above are the same as the first
nine exact values corresponding to elastic fixing with the
relative stiffness factors of 1, 2, 3, 4. This means that the
unknown rod fastening inaccessible to direct observation has been
correctly determined.

\end{document}